\documentclass[12pt]{amsart}
\usepackage{amsfonts,amscd,amsmath,amssymb,amsthm,mathrsfs}
\usepackage{array,latexsym,color,float,enumerate}
\usepackage{graphicx,epsfig}
\usepackage[top=2.7cm, bottom=2.7cm, left=3.2cm, right=3.2cm]{geometry}
\usepackage[all]{xy}

\makeindex

\usepackage{marginnote}

\theoremstyle{plain}
\newtheorem{theorem}{Theorem}[section]
\newtheorem{lemma}[theorem]{Lemma}
\newtheorem{proposition}[theorem]{Proposition}
\newtheorem{corrolary}[theorem]{Corollary}
\theoremstyle{definition}
\newtheorem{definition}[theorem]{Definition}
\newtheorem{example}[theorem]{Example}
\newtheorem{noname}[theorem]{}
\newtheorem{remark}[theorem]{Remark}
\newtheorem{construction}[theorem]{Construction}
\newtheorem{notation}[theorem]{Notation}
\theoremstyle{remark}
\newtheorem*{smallremark}{Remark}

\newtheorem{case}{Case} \makeatletter \@addtoreset{case}{theorem}\makeatother
\newtheorem{claim}{Claim} \makeatletter \@addtoreset{claim}{theorem}\makeatother

\newcommand{\bthm}{\begin{theorem}}
\newcommand{\bprop}{\begin{proposition}}
\newcommand{\blem}{\begin{lemma}}
\newcommand{\bcor}{\begin{corrolary}}
\newcommand{\brem}{\begin{remark}}
\newcommand{\bdfn}{\begin{definition}}
\newcommand{\bitem}{\begin{itemize}}
\newcommand{\bex}{\begin{example}}
\newcommand{\bno}{\begin{noname}}
\newcommand{\bsrem}{\begin{smallremark}}
\newcommand{\bnot}{\begin{notation}}
\newcommand{\bcon}{\begin{construction}}
\newcommand{\bca}{\begin{case}}
\newcommand{\bcl}{\begin{claim}}

\newcommand{\ecl}{\end{claim}}
\newcommand{\eca}{\end{case}}
\newcommand{\econ}{\end{construction}}
\newcommand{\enot}{\end{notation}}
\newcommand{\esrem}{\end{smallremark}}
\newcommand{\eno}{\end{noname}}
\newcommand{\eex}{\end{example}}
\newcommand{\eitem}{\end{itemize}}
\newcommand{\ethm}{\end{theorem}}
\newcommand{\eprop}{\end{proposition}}
\newcommand{\elem}{\end{lemma}}
\newcommand{\ecor}{\end{corrolary}}
\newcommand{\erem}{\end{remark}}
\newcommand{\edfn}{\end{definition}}
\newcommand{\benum}{\begin{enumerate}}
\newcommand{\eenum}{\end{enumerate}}
\newcommand{\un}{\underline}
\newcommand{\ov}{\overline}

\newcommand{\wt}{\widetilde}
\newcommand{\cal}[1]{\mathcal{#1}}

\newcommand{\ds}{\displaystyle}

\def\8{\infty}
\def\.{\cdot}
\def\PP{\mathbb{P}}

\def\C{\mathbb{C}}
\def\Z{\mathbb{Z}}
\def\N{\mathbb{N}}
\def\Q{\mathbb{Q}}

\def\E{\widehat{E}}
\def\ovk{\overline\kappa}

\def\xra{\xrightarrow}
\def\raa{\xrightarrow{\ \ }}

\def\:{\colon}

\def\bsk{\bigskip}

\def\Bk{\operatorname{Bk}}
\def\Proj{\operatorname{Proj}}

\def\Sing{\operatorname{Sing}}

\def\Ker{\operatorname{Ker}}
\def\im{\operatorname{Im}}

\def\Supp{\operatorname{Supp}}

\def\dim{\operatorname{dim}}
\def\Exc{\operatorname{Exc}}

\usepackage[hypertex,linkbordercolor={0 0 1}]{hyperref} 

\begin{document}
\title[Singular $\Q$-homology planes I]{Classification of singular $\Q$-homology planes. I.~Structure and singularities.}
\author[Karol Palka]{Karol Palka}
\address{Karol Palka: Institute of Mathematics, University of Warsaw, ul. Banacha 2, 02-097 Warsaw, Poland}
\address{Institute of Mathematics, Polish Academy of Sciences, ul. \'{S}niadeckich 8, 00-956 Warsaw, Poland}
\thanks{The author was supported by Polish Grant NCN N N201 608640}
\email{palka@impan.pl}

\begin{abstract} A $\Q$-homology plane is a normal complex algebraic surface having trivial rational homology. We obtain a structure theorem for $\Q$-homology planes with smooth locus of non-general type. We show that if a $\Q$-homology plane contains a non-quotient singularity then it is a quotient of an affine cone over a projective curve by an action of a finite group respecting the set of lines through the vertex. In particular, it is contractible, has negative Kodaira dimension and only one singular point. We describe minimal normal completions of such planes.    \end{abstract}

\maketitle

We work with algebraic varieties defined over $\C$.

\section{Main results}

A $\Q$-homology plane is a normal surface with Betti numbers of $\C^2$, i.e. with $H_i(-,\Q)=0$ for $i>0$. As for every open surface, one of its basic invariants is the logarithmic Kodaira dimension (see \cite{Iitaka}), which takes values in $\{-\8,0,1,2\}$. Smooth $\Q$-homology planes of non-general type, i.e. having the Kodaira dimension smaller than two, have been classified, see \cite[\S3.4]{Miyan-OpenSurf} for summary and for what is known in the case of general type. In this and in the forthcoming paper (see \cite{Palka-classification2}) we obtain a classification of singular $\Q$-homology planes with smooth locus of non-general type. A lot of attention has been given to understand these surfaces in special cases (see \cite{MiSu-hPlanes, GM-Affine-lines, PS-rationality, Daigle-Russell-logQhp, KR-ContrSurf}), let us mention explicitly at least the role of the contractible ones in proving the linearizability of $\C^*$-actions on $\C^3$ (see \cite{KR-C*onC3}). To our knowledge, in the available literature on this subject it is always assumed that the planes are \emph{logarithmic}, by what is meant that each singular point is analytically of type $\C^2/G$ for some finite subgroup $G<GL(2,\C)$ (singularity of 'quotient type'). This is a strong assumption, in particular it implies rationality of the surface (\cite{GPS-logQhp_rational}), and one of our goals was to avoid it.

Recall that a $\Q$-homology plane is \emph{exceptional} if and only if it has smooth locus of non-general type, which is neither $\C^1$- nor $\C^*$-ruled. There exist exactly three exceptional smooth $\Q$-homology planes (see \cite[8.64]{Fujita}). In \cite{Palka-k(S_0)=0} we described two singular ones. We obtain the following structure theorem, part (2) strengthening 1.4 loc. cit.

\begin{theorem}\ \label{thm:main result1} \benum[(1)]

\item Singular $\Q$-homology planes are affine and birationally ruled.

\item A singular $\Q$-homology plane with smooth locus of non-general type satisfies one of the following. \benum[(a)]

    \item It is logarithmic and $\C^1$- or $\C^*$-ruled.

    \item It is either non-logarithmic or isomorphic to $\C^2/G$ for some finite small non-cyclic subgroup $G<GL(2,\C)$. Its smooth locus has a $\C^*$-ruling which does not extend to a ruling of the $\Q$-homology plane.

    \item It is isomorphic to one of two exceptional singular $\Q$-homology planes. These planes have Kodaira dimension zero and the Kodaira dimension of the smooth loci zero. They are quotients of smooth exceptional $\Q$-homology planes and contain unique singular points, which are cyclic singularities of Dynkin type $A_1$ and $A_2$ respectively. \eenum

\eenum\end{theorem}

Logarithmic $\Q$-homology planes which admit a $\C^1$- or a $\C^*$-ruling have been studied in \cite{MiSu-hPlanes}, in particular singular fibers of the rulings have been described. In the forthcoming, second part of this paper we will give a detailed classification of these planes. As we see from the above, non-logarithmic $\Q$-homology planes are of type (b). We obtain the following more detailed description.

\begin{theorem}\label{thm:main result2} A singular $\Q$-homology plane which contains a non-quotient singularity (in particular each which is a non-rational surface) is isomorphic to a quotient of an affine cone over a smooth projective curve by an action of a finite group acting freely off the vertex of the cone and respecting the set of lines through the vertex. Moreover: \benum[(i)]

\item it is contractible and has a unique singular point,

\item its smooth locus has a unique $\C^*$-ruling and the ruling does not extend to a ruling of the plane,

\item it has negative Kodaira dimension and the Kodaira dimension of its smooth locus equals $0$ or $1$,

\item its minimal normal completion is unique up to isomorphism and the boundary is a rational tree with a unique branching component. \eenum\end{theorem}

We show also that for $\Q$-homology planes as above the singularity may be, but not need to be, rational in the sense of Artin (cf. \ref{ex4:nonextendable-singularities}). The following result is of independent interest (for the proof see \ref{prop3:general is logarithmic}).

\begin{proposition}\label{prop:main result} If a singular $\Q$-homology plane has smooth locus of general type then it has a unique singular point and the point is a quotient singularity. If a $\Q$-homology plane has more than one singular point then either it is affine-ruled and the singularities are cyclic or it has exactly two singular points, both of Dynkin type $A_1$. \end{proposition}

Our methods rely on the theory of open algebraic surfaces, for which \cite{Miyan-OpenSurf} is a basic reference. We now give a more detailed overview of the paper. We denote a singular $\Q$-homology plane by $S'$ and its smooth locus by $S_0$. From preliminaries in section \ref{sec:preliminaries} the result \ref{cor2:Surf/Div} is worth mentioning, as it gives a criterion for contractibility of divisors on complete surfaces. In section \ref{sec:topology} we prove basic topological and geometric results, whose simpler versions for logarithmic $\Q$-homology planes were known before, see \cite{MiSu-hPlanes}.  In absence of restriction on the type of singularities arguments get more complicated. Once we prove that the Neron-Severi group of the smooth locus is torsion, we apply Fujita's argument to show that the affiness of $S'$ is a consequence of $\Q$-acyclicity (cf. \ref{cor3:homology}). It was proved in \cite{PS-rationality} that a logarithmic singular $\Q$-homology plane is rational. We complete this result by showing in \ref{prop3:generalities on ovS-D-E} that in general $S'$ is birationally ruled.

If $S_0$ is of non-general type and admits no $\C^1$- and no $\C^*$-rulings then the plane $S'$ is called \emph{exceptional}. By general structure theorems for open surfaces (cf. \cite[2.2.1, 2.5.1, 2.6.1]{Miyan-OpenSurf}) exceptional $\Q$-homology planes have $\ovk(S_0)=0$. Under the assumption that singularities are topologically rational we have proved in \cite{Palka-k(S_0)=0} that there are exactly two such surfaces up to isomorphism. Here we show that the mentioned assumption can be omitted.

We begin section \ref{sec:nonlogarithmic S'} by proving that if $S'$ is non-logarithmic then it is of very special kind, namely there is a unique $\C^*$-ruling of $S_0$ and it does not extend to a ruling of $S'$. This implies $\ovk(S')=-\8$ (cf. \ref{thm3:nonlogarithmic are special}). We analyze the $\C^*$-ruling and classify all non-logarithmic $\Q$-homology planes in \ref{thm5:nonextendable}. To reconstruct them we use the contractibility criterion \ref{cor2:Surf/Div}. We analyze the Kodaira dimension of the smooth locus and the singularity in \ref{cor5:nonextendable-singularities} and \ref{ex4:nonextendable-singularities}. Finally we argue that a non-logarithmic $S'$ admits a $\C^*$-action with a unique fixed point, which leads to an isomorphism with a quotient of an affine cone as in \ref{thm:main result2}.

\bsk\textsl{\textsf{Acknowledgements.}} The paper contains results obtained during the graduate studies of the author at the University of Warsaw and their improvements obtained during his stay at the Polish Academy of Sciences and McGill University. The author thanks his thesis advisor prof. M. Koras for numerous discussions and for reading preliminary versions of the paper. He also thanks prof. P. Russell for useful comments.

\tableofcontents

\section{Preliminaries}\label{sec:preliminaries}

\subsection{Divisors and pairs}\label{ssec:div and pairs}

We mostly follow the notational conventions and terminology of \cite{Miyan-OpenSurf}, \cite{Fujita} and also of \cite[\S1, \S2]{Palka-k(S_0)=0}. Let $T=\sum t_iT_i$ be an snc-divisor on a smooth complete surface (hence projective by the theorem of Zariski) with distinct irreducible components $T_i$. We write $\un T=\sum T_i$ for a reduced divisor with the same support as $T$ and denote the branching number of $T_i$ by $\beta_T(T_i)=\un T\cdot (\un T-T_i)$. A component $T_i\subseteq T$ is \emph{branching} if $\beta_T(T_i)\geq 3$. If $T$ contains a branching component then it is \emph{branched}. The determinant of $-Q(T)$, where $Q(T)$ is the intersection matrix of $T$, is denoted by $d(T)$, $d(0)=1$ by definition. Considering $T$ as a topological subspace of a complex surface with its Euclidean topology it is easy to check that if $\Supp T$ is connected then there is a homotopy equivalence $$T\underset{htp}{\approx} \bigvee_{i=1}^n T_i\vee\ |DG(T)|,$$ where $|DG(T)|$ is a geometric realization of a dual graph of $T$. In particular, $b_1(T)=\ds\sum_{i=1}^n b_1(T_i)+b_1(|DG(T)|)$.

If $T$ is a chain (i.e. it is reduced and its dual graph is linear) then writing it as a sum of irreducible components $T=T_1+\ldots+T_n$ we always assume that $T_i\cdot T_{i+1}=1$ for $1\leq i\leq n-1$. If $T$ is a chain and some tip (a component with $\beta\leq 1$), say $T_1$, is fixed to be the first one then we distinguish between $T$ and $T^t=T_n+\ldots+T_1$. We write $T=[-T_1^2,\ldots,-T_n^2]$ in case $\un T$ is a rational chain. If $T$ is a rational chain with $T_i^2\leq -2$ for each $i$ we say that $T$ is \emph{admissible}. Let $D$ be some fixed reduced snc-divisor which is not an admissible chain. A rational chain $T\subseteq D$ not containing branching components of $D$ and containing one of its tips is a \emph{twig} of $D$. In this situation we always assume that the tip of $D$ is the first component of $T$. For any admissible (ordered) chain we define $$e(T)=\frac{d(T-T_1)}{d(T)}\text{\ and\ }\wt e(T)=e(T^t).$$ Now $e(D)$ and $\wt e(D)$ are defined as the sums of respective numbers computed for all maximal admissible twigs of $D$. (A an admissible twig of $D$ is maximal if it is not contained in another admissible twig of $D$ with a bigger number of components.)

An \emph{snc-pair} $(X,D)$ consists of a complete surface $X$ and a reduced snc-divisor $D$ contained in the smooth part of $X$. We write $X-D$ for $X\setminus D$ in this case. The pair is a \emph{normal pair} (\emph{smooth pair}) if $X$ is normal (resp. smooth). If $X$ is a normal surface then an embedding $\iota\:X\to \ov X$, where $(\ov X,\ov X\setminus X)$ is a normal pair, is called a \emph{normal completion} of $X$. If $X$ is smooth then $\ov X$ is smooth and $(\ov X,D,\iota)$ is called a \emph{smooth completion} of $X$. We often identify $X$ with $\ov X-D$ by $\iota$ and neglect $\iota$ in the notation. A morphism of two completions $\iota_j:X\to \ov X_j$, $j=1,2$ is a morphism $f\:\ov X_1\to \ov X_2$, such that $\iota_2=f\circ\iota_1$.

Let $\pi\:(X,D)\to (X',D')$ be a birational morphism of normal pairs. We put $\pi^{-1}D'=\un {\pi^*D'}$, i.e. $\pi^{-1}D'$ is the reduced total transform of $D'$. Assume $\pi^{-1}D'=D$. If $\pi$ is a blowup then we call it \emph{subdivisional} (\emph{sprouting}) for $D'$ if its center belongs to two (one) components of $D'$. In general we say that $\pi$ is \emph{subdivisional} for $D'$ (and for $D$) if for any component $T$ of $D'$ we have $\beta_{D'}(T)=\beta_D(\pi^{-1}T)$.

The exceptional locus of a birational morphism between two surfaces $\eta:X\to X'$, denoted by $\Exc(\eta)$, is defined as the locus of points in $X$ for which $\eta$ is not a local isomorphism. The canonical divisor of a complete surface $X$ is denoted by $K_X$ and the numerical equivalence of divisors by $\equiv$. For a divisor $D$ the arithmetic genus of $D$ is $p_a(D)=\frac{1}{2}D\cdot(K_X+D)+1$.

A $b$-curve is a smooth rational curve with self-intersection $b$. A divisor is snc-minimal if all its $(-1)$-curves are branching.

\bdfn\label{def2:connected morphism} A birational morphism of surfaces $\pi:X\to X'$ is a \emph{connected modification} if it is proper, $\pi(\Exc(\pi))$ is a smooth point on $X'$ and $\Exc(\pi)$ contains a unique $(-1)$-curve. In case $\pi$ is a morphism of pairs $\pi\:(X,D)\to (X',D')$, such that $\pi^{-1}(D')=D$ and $\pi(\Exc(\pi))\in D'$, then we call it a \emph{connected modification over $D'$}.\edfn

Note that since for a connected modification the exceptional locus contains a unique $(-1)$-curve, the modification can be decomposed into a sequence of blowdowns $\sigma_n\circ\ldots\circ\sigma_1$, such that for $i\leq n-1$ the center of $\sigma_i$ belongs to the exceptional divisor of $\sigma_{i+1}$. A sequence of blowdowns (and its reversing sequence of blowups) whose composition is a connected modification will be called a \emph{connected sequence of blowdowns (blowups)}.

\blem\label{lem2:Eff-NegDef is Eff} Let $A$ and $B$ be $\Q$-divisors on a smooth complete surface, such that $Q(B)$ is negative definite and $A\cdot B_i\leq 0$ for each irreducible component $B_i$ of $B$. Denote the integral part of a $\Q$-divisor by $[\ ]$.\benum[(i)]

\item If $A+B$ is effective then $A$ is effective.

\item If $n\in \N$ and $n(A+B)$ is a $\Z$-divisor then $h^0(n(A+B))=h^0([nA])$.\eenum\elem

\begin{proof} (i) We can assume that $A$ and $B$ are $\Z$-divisors and $B$ is effective and nonzero. Write $B=\sum b_iB_i$, where $B_i$ are distinct irreducible components of $B$. Choose $b_i'\in \N$, such that the sum $\sum b_i'$ is the smallest possible among divisors $\sum b_i'B_i$, such that $A+\sum b_i'B_i$ is effective. If $b_i'>0$ for some $i$ then $$(A+\sum b_i'B_i)\cdot(\sum b_i'B_i)\leq(\sum b_i'B_i)^2<0$$ by the assumptions. Hence $\Supp(A+\sum b_i'B_i)$ contains some $B_i$, a contradiction with the definition of $b_i'$. Thus $A$ is effective.

(ii) Let $\{R\}$ denote the fractional part of a $\Q$-divisor $R$, i.e. $\{R\}=R-[R]$. Let $T$ be some effective divisor, such that $n(A+B)\sim T$. Then $nA\sim T-nB$ as $\Q$-divisors. Since $T-nB$ is effective by (i), the coefficient of each irreducible component of $[T-nB]$ is bounded below by the coefficient of the same component in $-\{T-nB\}$. Since $[T-nB]$ is a $\Z$-divisor and the coefficients of components in $\{T-nB\}$ are fractional and positive, $[T-nB]$ is effective. Moreover, $\{T-nB\}-\{nA\}$ being a $\Z$-divisor is equal to $0$, so the rational function giving the equivalence of $n(A+B)$ and $T$ gives an equivalence of $[T-nB]$ and $[nA]$. \end{proof}

\subsection{Singularities and contractible divisors}\label{ssec:sing}

Let $\E$ be the reduced exceptional divisor of the (unique) minimal good (i.e. such that $\E$ is an snc-divisor) resolution  of a singular point on a normal surface $X$. Then $\E$ is connected and $Q(\E)$ is negative definite. Recall that a point $q\in X$ is of \emph{quotient type} if there exists an analytical neighborhood $N\subseteq X$ of $q$ and a small (i.e. not containing any pseudo-reflections) finite subgroup $G$ of $GL(2,\C)$, such that $(N,q)$ is analytically isomorphic to $(\wt{N}/G,0)$ for some ball $\wt{N}$ around $0\in \C^2$. Then $G=\pi_1(N\setminus\{q\})$. Note that by a result of Tsunoda (\cite{Tsunoda-structure_of_open_surf}) for normal surfaces quotient singularities are the same as log-terminal singularities. For a singular point $q\in X$ of quotient type it is known (\cite{Brieskorn}) that $G$ is cyclic if and only if $\E$ is an admissible chain and that $G$ is non-cyclic if and only if it is non-abelian if and only if $\E$ is an admissible fork (rational snc-minimal fork with three twigs and with negative definite intersection matrix, cf. \cite[2.3.5]{Miyan-OpenSurf}), in each case $d(\E)=|G/[G,G]|$. In case $\E$ is a fork, we will say that $\E$ is \emph{of type $(d_1,d_2,d_3)$} if the maximal twigs of $\E$ have $d(\ )$ equal to $d_1, d_2, d_3$. Quotient singularities are \emph{rational}, as the first direct image of the structure sheaf of their resolutions vanishes. It follows from \cite[1]{Artin-Rational} that for a rational singularity $\E$ is a rational tree, hence rational singularities are \emph{topologically rational}, which by definition means that $b_1(\E)=0$. This notion is a bit stronger than the \emph{quasirationality} in the sense of Abhyankar (cf. \cite{Ab}), for which only the rationality of components of $\E$ is required.

\bex\label{ex1:Pham-Brieskorn} Let $V\subseteq\C^3$ be given by $x^2+y^3+z^7=0$. Then the blowup of $V$ in $0$ has an exceptional line contained in the singular locus, hence is not normal. Since the blowup of a normal surface with rational singularity remains normal by \cite[8.1]{Lipman}, $0\in V$ is not a rational singularity. On the other hand, it is topologically rational.

More generally, let $V(p_1,p_2,p_3)\subseteq \C^3$ be a Pham-Brieskorn surface given by the equation $x_1^{p_1}+x_2^{p_2}+x_3^{p_3}=0$, where $p_1,p_2,p_3\geq 2$. This surface is contractible (note it has a $\C^*$-action with the singularity as the unique fixed point) and it is known that $0\in V(p_1,p_2,p_3)$ is a topologically rational singularity if and only if one of $p_1,p_2,p_3$ is coprime with two others or $\frac{1}{2}p_1,\frac{1}{2}p_2,\frac{1}{2}p_3$ are integers coprime in pairs. (In \cite{Orevkov} and \cite[0.1]{Zaidenberg} the above is stated as a condition for quasirationality, but in both cases the graph of the resolution is a tree by looking at the proof or by using \cite{Orlik-Wagreich}). On the other hand, the rationality of $0\in V(p_1,p_2,p_3)$ is by \cite[2.21]{Zaidenberg} equivalent to each of the following conditions: \benum[(i)]
\item $\sum_{i=1}^3\frac{1}{p_i}>1$, \item $0\in V(p_1,p_2,p_3)$ is of quotient type, \item $\ovk(V\setminus\{0\})=-\8$. \eenum\eex

We have the following corollary from the Nakai criterion.

\blem\label{lem2:contraction lemma} Let $A$ and $B$ be effective snc-divisors on a smooth complete surface $X$ having disjoint supports. If for every irreducible curve $C$ on $X$ either $C\subseteq B$ or $A\cdot C>0$ then for sufficiently large and sufficiently divisible $n$ one has:\benum[(i)]

\item $|nA|$ has no base points,

\item $\varphi_{|nA|}$ is birational and contracts exactly the curves in $B$,

\item $\varphi_{|nA|}(X)$ is normal, projective and isomorphic to $$\Proj \underset{n\geq0}{\bigoplus} H^0(\cal O_X(nA)).$$ \eenum\elem

\begin{proof} (i) Repeating part of the proof of Nakai's criterion (cf. \cite[V.1.10]{Har}) we get that $\cal O(nA)$ is generated by global sections for $n\gg 0$. For (ii) and (iii) see for example \cite[2.3, 2.4]{Reid}. See also \cite[3.4]{Schroer-contractible_curves} for contractibility criterion for normal surfaces not involving effectiveness. \end{proof}

\bdfn Let $(\ov X,D)$ be a smooth completion of a smooth surface $X$ and let $NS(\ov X)$ be the Neron-Severi group of $\ov X$ consisting of numerical equivalence classes of divisors. The Neron-Severi group $NS(X)$ of $X$ is defined as the cokernel of the natural map $\Z[D]\to NS(\ov X)$, where $\Z[D]$ is a free abelian group generated by irreducible components of $D$. We denote $NS(X)\otimes\Q$ by $NS_\Q(X)$.\edfn

\bsrem The above definition does not depend on a smooth completion of $X$ (cf. \cite[1.19]{Fujita}). Contrary to the case when $X$ is complete, in general $NS(X)$ can have torsion.\esrem

\bcor\label{cor2:Surf/Div} Let $A$ and $B$ be effective snc-divisors on a smooth complete surface $X$ having disjoint supports. Assume that $A$ is connected, $Q(B)$ is negative definite and $NS_\Q(X-A-B)=0$. Then there exists a normal affine surface $Y$ and a morphism $\zeta\:X-A\to Y$ contracting connected components of $B$, such that $\zeta\:X-A-B\to Y\setminus\zeta(B)$ is an isomorphism.\ecor

\begin{proof}
Smooth complete surface is projective by the theorem of Zariski. Since $NS_\Q(X-A-B)=0$, there exists a divisor $H=H_A+H_B$ with $H_A\subseteq A$ and $H_B\subseteq B$, which is numerically equivalent to an ample divisor on $X$. Then $H$ is ample, because ampleness is a numerical property by Nakai's criterion. To use \ref{lem2:contraction lemma} we need to show that there exists an effective divisor $F$, such that $\Supp F=\Supp A$ and $F\cdot C>0$ for all irreducible curves $C\nsubseteq B$. To deal with curves $C\subseteq A$ we use Fujita's argument (\cite[2.4]{Fujita}). Let $\mathcal{U}$ consist of all effective divisors $T$, such that $T\subseteq A$ and $T\cdot T_i>0$ for any prime component $T_i$ of $T$. Writing $H_A=H_+-H_-$, where $H_+$, $H_-$ are effective and have no common component, we see that $\mathcal{U}$ is nonempty because $H_+\in\mathcal{U}$. Suppose $F$ is an element of $\mathcal{U}$ with maximal number of components. For an irreducible curve $C\nsubseteq F$ satisfying $C\cdot F>0$ one would get $tF+C\in\mathcal{U}$ for $t>max(0,-C^2)$, hence $\Supp F=\Supp A$ by the connectedness of $A$.

Suppose an irreducible curve $C\nsubseteq B$ satisfies $C\cdot F=0$. Since $F\in \mathcal{U}$, we have $C\nsubseteq F$. We can choose some reduced divisor $F'\subseteq F$, such that irreducible components of $F'+B$ give a basis of $NS_\Q(X)$. Let us write $C\equiv \sum_i\alpha_i F_i+B^+-B^-$, where $F_i\subseteq F'$, the divisors $B^+,B^-\subseteq B$ are effective and have no common component. For each $j$ we have $C\cdot F_j=0$, so $(\sum_i\alpha_i F_i)\cdot F_j=C\cdot F_j=0$, hence $\sum_i\alpha_i F_i=0$ because $d(F')\neq0$. We have $$(B^+)^2=B^+\cdot C+B^+\cdot B^-\geq 0,$$ so $B^+=0$. Thus the divisor $C+B^-$ is nonzero, effective and numerically trivial, a contradiction. Let $\zeta=\varphi_{|nF|}$ for $n$ as in lemma \ref{lem2:contraction lemma}. Then $\zeta$ contracts connected components of $B$. We have also $nF=\zeta^*H$, where $H$ is a very ample divisor on $\zeta(X)$, which implies that $\zeta(X-A)$ is affine. \end{proof}

\bsrem Note that any divisor with negative definite intersection matrix can be contracted in the analytical category by the theorem of Grauert (cf. \cite{Grauert}). However, in general it is a more subtle problem if this can be done in the algebraic category (see \cite{Artin-Rational} for results concerning rational singularities). \esrem

\subsection{Minimal models}

Let us give a brief sketch of the notion of minimality for open surfaces for unfamiliar readers. By the Castelnuovo criterion a smooth projective surface $X$ is minimal if and only if there is no irreducible curve $L$ on $X$ for which $K_X\cdot L<0$ and $L^2<0$, which is equivalent to $L$ being a $(-1)$-curve. Similarly, we can say that a smooth pair $(X,D)$ is \emph{relatively minimal} if and only if there is no irreducible curve $L$ on $X$ for which $(K_X+D)\cdot L<0$ and $L^2<0$. In case $L\nsubseteq D$ this implies that $L$ is a $(-1)$-curve intersecting $D$ in at most one point and transversally. However, if $L\subseteq D$ then the conditions are equivalent to $L^2<0$ and $$\beta_D(L)=L\cdot (D-L)<2(1-p_a(L))$$ and hence to $L$ being a smooth rational curve with negative self-intersection and branching number $\beta_D(L)<2$. Contraction of such an $L$ immediately leads out of the category of smooth pairs, as in particular any tip of any admissible maximal twig of $D$ would have to be contracted. Thus one repeats the definition of a relatively minimal pair for pairs $(\ov X,\ov D)$ consisting of a normal projective surface and reduced Weil divisor (cf. \cite[2.4.3]{Miyan-OpenSurf}). Then a relatively minimal model of a given pair (which can be singular and not unique) is obtained by successive contractions of curves satisfying the above conditions. To go back to the smooth category one can translate the conditions for $(\ov X,\ov D)$ to be relatively minimal in terms of the properties of its minimal resolution. This leads to the notion of an \emph{almost minimal pair}, which we recall now for the convenience of the reader (cf. \cite[2.3.11]{Miyan-OpenSurf}).

First, for any smooth pair $(X,D)$ we define the \emph{bark of $D$}. For non-connected $D$ bark is a sum of barks of its connected components, so we will assume $D$ is connected. If $D$ is an snc-minimal resolution of a quotient singularity (i.e. $D$ is an admissible chain or an admissible fork) then we define $\Bk D$ as a unique $\Q$-divisor with $\Supp \Bk D\subseteq D$, such that  $$(K_X+D-\Bk D)\cdot D_i=0\text{\ for each component\ } D_i\subseteq D.$$ In other case let $T_1,\ldots, T_s$ be all the maximal admissible twigs of $D$. (If $\ovk(X-D)\geq 0$ and $D$ is snc-minimal then all rational maximal twigs of $D$ are admissible, cf. \cite[6.13]{Fujita}). In this case we define $\Bk D$ as a unique $\Q$-divisor with $\Supp \Bk D\subseteq \bigcup T_j$, such that $$(K_X+D-\Bk D)\cdot D_i=0\text{\ for each component\ } D_i\subseteq \bigcup_{j=1}^s T_j.$$ The definition implies that $\Bk D$ is an effective $\Q$-divisor with negative definite intersection matrix and its components can be contracted to quotient singular points. In fact all components of $\Bk D$ in its irreducible decomposition have coefficients at smaller than $1$, unless $D$ is an admissible chain or fork consisting of $(-2)$-curves. Thus if $D$ is not such a $(-2)$-chain or a $(-2)$-fork then $D^\#=D-\Bk D$ is an effective divisor with $\Supp D^\#=\Supp D$.

We now say that a smooth pair $(X,D)$ is \emph{almost minimal} if for each curve $L$ on $X$ either $(K_X+D^\#)\cdot L\geq 0$ or $(K_X+D^\#)\cdot L<0$ but the intersection matrix of $\Bk D+L$ is not negative definite. Consequently, an almost minimal model of a given pair $(X,D)$ can be obtained by successive contractions of curves $L$
for which $(K_X+D^\#)\cdot L<0$ and $\Bk D+L$ is negative definite. These are the non-branching $(-1)$-curves in $D$ and $(-1)$-curves $L\nsubseteq D$ for which $D^\#\cdot L<1$ and $\Bk D+L$ is negative definite. Minimalization does not change the logarithmic Kodaira dimension. One shows that $(X,D)$ is almost minimal if and only if after taking the contraction $\epsilon:(X,D)\to (\ov X,\ov D)$ of connected components of $\Bk D$ to singular points the pair $(\ov X,\ov D)$ is relatively minimal. Moreover, if $(X,D)$ is almost minimal and $\ovk(X-D)\geq 0$ then $K_X+D^\#$ and $\Bk D$ are the numerically effective (nef) and negative definite parts of the Zariski decomposition of $K_X+D$. The reader can find a more detailed description of the process of minimalization in loc. cit. We will need in particular the following fact.

\brem\label{rem:almost minimalization} Let $(X,D)$ be a smooth pair which is not almost minimal, but for which $D$ is snc-minimal. Let $L\subseteq X$ be a curve witnessing the non-almost-minimality, i.e. $L$ is a $(-1)$-curve not contained in $D$, such that $D^\#\cdot L<1$ and $\Bk D+L$ is negative definite. Then $L$ meets $D$ transversally, in at most two points, each connected component of $D$ at most once. Moreover, if $L$ meets $D$ in two points then one of the connected components is an admissible chain and both points of intersection belong to $\Supp \Bk D$. \erem

\subsection{Rational rulings}\label{ssec:rulings}

By a \emph{rational ruling} of a normal surface we mean a surjective morphism of this surface onto a smooth curve, for which a general fiber is a rational curve. If its general fiber is isomorphic to $\PP^1$ it is called a $\PP^1$-ruling.

\bdfn\label{def2:morphism completion} If $p_0:X_0\to B_0$ is a rational ruling of a normal surface then by a \emph{completion of $p_0$} we mean a triple $(X,D,p)$, where $(X,D)$ is a normal completion of $X_0$ and $p\:X\to B$ is an extension of $p_0$ to a $\PP^1$-ruling with $B$ being a smooth completion of $B_0$. We say that $p$ is a \emph{minimal completion of $p_0$} if $D$ is \emph{$p$-minimal}, i.e. if $p$ does not dominate any other completion of $p_0$.\edfn

Note that if $D$ and $p$ are as above then $D$ is $p$-minimal if and only if each non-branching $(-1)$-curve contained in $D$ is horizontal.

For any rational ruling $p_0$ as above there is a completion $(X,D,p)$. Let $f$ be a general fiber of $p$. We call $p_0$ a $\C^1$-ruling (a $\C^{(n*)}$-ruling) if $f\cdot D=1$ (if $f\cdot D=n+1$). Any fiber of a $\PP^1$-ruling has vanishing arithmetic genus and self-intersection. The following well-known lemma shows that these conditions are also sufficient.

\blem\label{lem1:0-curve gives a ruling} Let $F$ be a connected snc-divisor on a smooth complete surface $X$. If $p_a(F)=F^2=0$ then there exists a $\PP^1$-ruling $p\:X\to B$ and a point $b\in B$ for which $p^*b=F$. \elem

\begin{proof} The proof given in \cite[V.4.3]{BHPV} after minor modifications works with the above assumptions. \end{proof}

For an irreducible vertical curve $J$ we denote its multiplicity in the fiber $F$ containing it by $\mu_F(J)$ (or $\mu(J)$ if $F$ is fixed).

If $F$ is a singular fiber of a $\PP^1$-ruling then, since $p_a(F)=0$, it is a rational tree with components of negative self-intersection. Its structure is well-known (see \cite[\S 4]{Fujita}). First of all, since $K_X\cdot F=-2$ (by the adjunction formula $F^2+K_X\cdot F=-2$), it contains a $(-1)$-curve and if the $(-1)$-curve is unique then its multiplicity is bigger than one. In fact each $(-1)$-curve of $F$ intersects at most two other components of $F$, so its contraction leads to a snc-fiber with a smaller number of components and, by induction, to a smooth fiber. In the process of contraction the total number of $(-1)$-curves in a fiber decreases, unless $F=[2,1,2]$.

The situation when $F$ has a unique $(-1)$-curve, say $C$, is of special interest. In this case $F$ is produced by a connected sequence of blowups from a smooth $0$-curve. Let $B_1,\ldots, B_n$ be the branching components of $F$ written in order in which they are created, put $B_{n+1}=C$. It is convenient to write $\un F$ as $\un F=T_1+T_2+ \ldots+T_{n+1}$, where the divisor $T_i$ is a reduced chain consisting of all components of $\un F-T_1-\ldots -T_{i-1}$ created not later than $B_i$. We call $T_i$ the \emph{i-th branch} of $F$. The proof of the following result is straightforward.

\blem\label{lem2:fiber with one exc comp} Let $F$ be a singular fiber of a $\PP^1$-ruling of a smooth complete surface. If $F$ contains a unique $(-1)$-curve $C$ then:\benum[(i)]

\item $\mu(C)>1$ and there are exactly two components of $F$ with multiplicity one. They are tips of the fiber and belong to the first branch,

\item if $\mu(C)=2$ then either $F=[2,1,2]$ or $C$ is a tip of $F$ and then $\un F-C=[2,2,2]$ or  $\un F-C$ is a $(-2)$-fork of type $(2,2,n)$,

\item if $F$ is branched then the connected component of $\un F-C$ not containing curves of multiplicity one is a chain (possibly empty).\eenum\elem

\bnot\label{not:rulings} Recall that having a fixed $\PP^1$-ruling of a smooth surface $X$ and a divisor $D$ we define  $$\Sigma_{X-D}=\underset{F\nsubseteq D}{\sum}(\sigma(F)-1),$$ where $\sigma(F)$ is the number of $(X-D)$ -components of a fiber $F$ (cf. \cite[4.16]{Fujita}). The horizontal part $D_h$ of $D$ is a divisor without vertical components, such that $D-D_h$ is vertical. The numbers $h$ and $\nu$ are defined respectively as $\#D_h$ and as the number of fibers contained in $D$. We will denote a general fiber by $f$. \enot

With the notation as above the following equation is satisfied (cf. loc. cit. and \cite[2.2]{Palka-k(S_0)=0} for a short proof): $$\Sigma_{X-D}=h+\nu+b_2(X)-b_2(D)-2.$$

\bdfn\label{def2:(un)twisted and columnar fiber} Let $(X,D,p)$ be a completion of a $\C^*$-ruling of a normal surface $X$. We say that the original ruling $p_0=p_{|X-D}$ is \emph{twisted} if $D_h$ is a 2-section. If $D_h$ consists of two sections we say that $p_0$ is \emph{untwisted}. A singular fiber $F$ of $p$ is \emph{columnar} if and only if it is a chain not containing singular points of $X$ and which can be written as $$\un F=A_n+\ldots+A_1+C+B_1+\ldots+B_m$$ with a unique $(-1)$-curve $C$, such that $D_h$ meets $F$ exactly in $A_n$ and $B_m$, in each once and transversally. The chains $A=A_1+\ldots+A_n$ and $B=B_1+\ldots+B_m$ are called \emph{adjoint chains}. \edfn

\bsrem By \cite[2.1.1]{KR-ContrSurf} and the fact that $d(A)$ and $d(A-A_1)$ are coprime we get easily that $e(A)+e(B)=1$ and $d(A)=d(B)=\mu_F(C)$. In fact we have also $\wt e(B)+\wt e(A)=1$ (see \cite[3.7]{Fujita}). \esrem

By abuse of language we call $p$ twisted or untwisted depending on the type of $p_0$. Twisted and untwisted $\C^*$-rulings are called respectively \emph{gyoza} (a Chinese dumpling) and \emph{sandwich} in \cite{Fujita}.

\section{\label{sec:topology}Topology and Singularities}

\subsection{Homology}\label{ssec:homology}

Let $S'$ be a \emph{singular $\Q$-homology plane}. Let $\epsilon \colon S \to S'$ be a good resolution and $(\ov S,D)$ a smooth completion of $S$. Denote the singular points of $S'$ by $p_1,\ldots,p_q$ and the smooth locus by $S_0$. We put $\E_i=\epsilon^{-1}(p_i)$ and assume that $\E =\E_1+\E_2+\ldots +\E_q$ is snc-minimal. Define $M$ as the boundary of the closure of $Tub(\E)$, where $Tub(\E)$ is a tubular neighborhood of $\E$. The construction of $Tub(\E)$ can be found in \cite{Mumford}. We may assume that $M$ is a disjoint sum of  $q$ closed oriented 3-manifolds. We write $H_i(X,A)$ for $H_i(X,A,\Q)$ and $b_i(X,A)$ for $\dim H_i(X,A)$.

Let us mention that the results we obtain below are generalizations of similar results obtained in the logarithmic case by Miyanishi and Sugie. However, restriction to quotient singularities is a strong assumption, which makes the considerations easier, even if at the end we prove that not so many non-logarithmic $\Q$-homology planes do exist.

\bprop\label{prop3:homology} Let $j_{\E}\:\E \to S$, $j_M\: M \to S_0$, $i_D\: D\to \ov S$ and $i_{D\cup\E}\:D\cup\E\to\ov S$ be the inclusion maps. One has:\benum[(i)]

\item $H_1(M,\Z)=H_1(\E,\Z)\oplus K$ for some finite group $K$ of order $d(\E)$,

\item $H_k(j_{\E})$ and $H_k(j_M)$ are isomorphisms for positive $k$,

\item $D$ is connected, $H_1(i_D)$ is an isomorphism and $b_1(D)=b_1(\E)$,

\item $H_2(i_{D\cup\E})$ is an isomorphism,

\item $H_k(S',\Z)=0$ for $k\neq 0,1$,

\item $\pi_1(\epsilon)\:\pi_1(S)\to \pi_1(S')$ is an epimorphism, it is an isomorphism if $b_1(\E)=0$,

\item if $b_1(\E)=0$ then $|d(D)|=|d(\E)|\cdot|H_1(S',\Z)|^2$. \eenum\eprop

\begin{proof}

(i) By \cite{Mumford} there is an exact sequence $$0\raa K \raa H_1(M,\Z)\xrightarrow{r} H_1(\E,\Z)\raa 0,$$ where $K$ is a finite group of order $d(\E)$ and $r$ is induced by the composition of embedding of $M$ into the closure of $Tub(\E)$ with retraction onto $\E$. Since $H_1(\E,\Z)$ is free abelian, it follows that $H_1(M,\Z)=H_1(\E,\Z)\oplus K$.

(ii) Let $k>0$. We look at the reduced homology exact sequence of the pair $(S,\E)$. The pairs $(S,\E)$ and $(S',\Sing S')$ are 'good CW-pairs' (see \cite[Thm 2.13]{Hat}), so for $k\neq 1$ we have $$H_k(S,\E)=H_k(S',Sing\ S')=0$$ and then $H_k(j_{\E})\colon H_k(\E) \to H_k(S)$ induced by $j_{\E}$ is an isomorphism for $k>1$. Now $$b_1(S,\E)=b_1(S',\Sing S')=b_0(\E)-1,$$ so $b_1(S)=b_1(\E)$ and $H_1(j_{\E})$ is also an isomorphism. Since $H_k(j_{\E})$ are epimorphisms, the Mayer-Vietories sequence for $S=S_0\cup Tub(\E)$ splits into exact sequences: $$0\raa H_k(M) \raa H_k(S_0)\oplus H_k(\E)\raa H_k(S)\raa 0.$$ Since $H_k(j_{\E})$ is injective, $H_k(j_M)$ is injective by exactness, so it is an isomorphism.

(iii)-(iv) By (ii) $b_3(S)=b_4(S)=0$, so the homology exact sequence of the pair $(\ov S,S)$ yields $H_4(\ov S,S)\cong H_4(\ov S)$, hence $H^0(D)=H_4(\ov S,S)=\Q$ by the Lefschetz duality (see \cite[7.2]{Dold}), which implies the connectedness of $D$. The components of $\E$ are numerically independent because $d(\E)\neq 0$, hence they are independent in $H_2(\ov S)$. This implies that the inclusion $i_{\E}\:\E\to \ov S$ induces a monomorphism on $H_2$. By (ii) we can write the exact sequence of the pair $(\ov S,S)$ as: $$\hdots\raa 0\raa H_3(\ov S) \raa H_3(\ov S,S) \raa H_2(\E) \raa H_2(\ov S)\to\hdots .$$ Now $H_2(i_{\E})$ is a monomorphism, so by the Lefschetz and Poincare duality $$b_1(D)=b_3(\ov S,S)=b_3(\ov S)=b_1(\ov S).$$ On the other hand $b_1(\ov S,D)=b_3(S)=0$, so $H_1(i_D)$ is an isomorphism.

Since $H_1(i_D)$ is an isomorphism, the homology exact sequence of the pair $(\ov S,D\cup\E)$ yields an exact sequence: $$0 \raa H_3(\ov S) \raa H_3(\ov S,D\cup\E) \xra{\ \delta\ } H_2(D\cup\E) \xra{\ \gamma\ } H_2(\ov S) \raa H_2(\ov S,D\cup\E) \raa H_1(\E)\raa 0.$$ We have $b_2(\ov S,D\cup\E)=b_2(S_0)=b_2(M)$ by (ii) and $b_2(M)=b_1(M)=b_1(\E)$ by (i), so $\gamma$ is an epimorphism. We need to prove that $b_1(D)=b_1(\E)$ and $\Ker \gamma =0$. Note that $b_2(\ov S)=b_2(D\cup\E)-\dim \im \delta$ and $$\dim \im \delta=b_3(\ov S,D\cup\E)-b_3(\ov S)=b_1(S_0)-b_1(\ov S)=b_1(\E)-b_1(D),$$ so $b_2(D\cup\E)-b_2(\ov S)=b_1(\E)-b_1(D)$. This implies that $b_1(D)=b_1(\E)$ if and only if $\Ker \gamma=0$.

If $b_1(\E)=0$ then $$b_3(\ov S,D\cup \E)=b_1(S_0)=b_1(\E)=0,$$ so $\gamma$ is a monomorphism. We can therefore assume that $\E$ is not a rational forest, in particular $S'$ is not logarithmic. Note that since $\gamma$ is an epimorphism, $S'$ is affine by \ref{cor2:Surf/Div}, so we can use \ref{prop3:general is logarithmic}(i) below to infer that $\ovk(S_0)\neq 2$. Suppose $\ovk(S_0)=1$, then $S_0$ is $\C^*$-ruled (cf. \cite[2.3]{Kawamata2}). Since modifications over $D+\E$ do not change $b_1(D)$ and $b_1(\E)$, we can assume that this ruling extends to $\ov S$. The divisor $D$ is not vertical, otherwise $Q(D+\E)$ would be semi-negative definite, which contradicts the Hodge index theorem. On the other hand, $\E$ is not vertical because is not a rational forest, so each of $D$ and $\E$ contains a unique section. Then $b_1(D)=b_1(\E)$, so we are done. We can now assume $\ovk(S_0)\leq 0$. Suppose $\ovk(S)=0$, then $\ovk(S_0)=0$. Put $F=D+\E-E_0$, where $E_0$ is a connected component of $\E$ with $b_1(E_0)\neq 0$. Let $(\wt S,\wt F+\wt E_0)$ be an almost minimal model of $(\ov S,D+\E)=(\ov S,F+E_0)$ with $\wt F$ and $\wt E_0$ being the direct images of $F$ and $E_0$. Each curve contracted in the process of minimalization intersects the image of $D$, because $S'$ is affine. Note that not all components of $D$ are contracted in this process, as $D$ is not negative definite by the Hodge index theorem. Moreover, by \ref{rem:almost minimalization} such a curve cannot intersect a connected component of $\E$ with nontrivial $b_1$. Thus the divisors $\wt F$ and $\wt E_0$ are disjoint, so $$K_{\wt S}+\wt F-\Bk \wt F+\wt E_0-\Bk \wt E_0\equiv (K_{\wt S}+\wt F+\wt E_0)^+\equiv 0.$$ Since by \ref{lem2:Eff-NegDef is Eff}(ii) $$h^0([n(K_{\wt S}+\wt F-\Bk \wt F)])=h^0(n(K_{\wt S}+\wt F))\geq h^0(n(K_S+F)),$$ we have $K_{\wt S}+\wt F-\Bk \wt F\geq _{\Q} 0$, so $\wt E_0=\Bk \wt E_0$, which contradicts $b_1(\wt E_0)=b_1(E_0)\neq 0$. We get $\ovk(S)=-\8$, so $S$ is $\C^1$-ruled by \cite[2.2.1]{Miyan-OpenSurf}. Consider an extension of this ruling to $\ov S$ and a divisor $$T=\sum_i d_i D_i+\sum_j e_jE_j$$ with distinct irreducible components $D_i\subseteq D$, $E_j\subseteq \E$, such that $T\equiv 0$. To finish the proof that $\Ker \gamma=0$ it is enough to show that $T=0$. Suppose $T\neq 0$. Using negative definiteness of $Q(\E)$ we see that each $e_j$ vanishes, otherwise $$0>(\sum_j e_jE_j)^2=T\cdot (\sum_j e_jE_j).$$ Intersecting $T$ with a fiber we see that the horizontal component of $D$ does not occur in the sum $T=\sum_j d_jD_j$, therefore $T$ is vertical. It follows that $\Supp T$ contains at least one fiber, otherwise $T^2<0$. However, then the equality $T\cdot \E=0$ implies that $\E$ is vertical, a contradiction with $b_1(\E)\neq 0$.

(v) Let $k\in\{3,4\}$. The groups $H_k(S',\Z)\cong H_k(S,\E,\Z)$ are torsion, so the exact sequence of the pair $(S,\E)$ gives $H_k(S,\E,\Z)\cong H_k(S,\Z)$. By the universal coefficient formula and Lefschetz duality $$H_k(S,\Z)\cong H^{k+1}(S,\Z)\cong H_{3-k}(\ov S,D,\Z)=0.$$ Vanishing of $H_2(S',\Z)$ is more subtle. The generalization of Andreotti-Frankel theorem to the singular case proved by Kar{\v{c}}jauskas says that an affine variety $X$ of complex dimension $n$ has the homotopy type of a $CW$-complex of real dimension not greater than $n$ (see \cite{GorMcP} for proofs and generalizations). In particular, $H_n(X,\Z)$ is a free abelian group. We showed in the proof of (iii)-(iv) that $S'$ is affine, so we get $H_2(S',\Z)=0$.

(vi) Choose points $y,x_1,\ldots,x_q\in S$, such that $y\in S_0$ and $x_i\in \E_i$. Let $\alpha_1,\ldots,\alpha_q$ be smooth paths in $S$ joining $y$ with $x_i$. We can choose $\alpha_i$ in such a way that they meet transversally in $y$, $\alpha_i\setminus \{y\}$ are disjoint, $R=\bigcup_i\alpha_i\setminus~\{x_i\}$ is contained in $S_0$ and meets $\E$ transversally. Let $N$ be a tubular neighborhood of $\E\cup R$ in $S$. Then $\epsilon(N)\subseteq S'$ is a contractible neighborhood of $\Sing S'$. Put $H=\pi_1(N\setminus (R\cup\E))$. Clearly, $\epsilon$ identifies $N\setminus (R\cup\E)$ with $\epsilon(N)\setminus(\Sing S'\cup\epsilon(R))$, so since $\pi_1(S_0\setminus R)\cong \pi_1(S_0)$, by van Kampen's theorem $$\pi_1(S)\cong \pi_1(S_0)\underset{H}{*}\pi_1(N)$$ and $$\pi_1(S')\cong \pi_1(S_0)\underset{H}{*}\{1\}.$$ We have $\pi_1(N)=\pi_1(\E_1)*\ldots*\pi_1(\E_q)$ and each $\pi_1(\E_i)$ is contained in the kernel of $\pi_1(\epsilon)$. If $b_1(\E)=0$ then $\E$ is a rational forest, so $\pi_1(N)=\{1\}$ and we get $\pi_1(S')\cong \pi_1(S)$.

(vii) Let $M_D=\partial Tub(D)$ be the boundary of a (closure of a) tubular neighborhood of $D$. We may assume that $M_D$ is a 3-manifold disjoint from $M$. By (iii) $b_1(D)=0$ and by (iv) the components of $D$ are independent in $H_2(\ov S)$, so $D$ is a rational tree with $d(D)\neq 0$. Then using the presentation given in \cite{Mumford} we get that $H_1(M_D)$ is a finite group of order $|d(D)|$. By Poincare duality $H_2(M_D,\Z)$ (and similarly $H_2(M,\Z)$) is trivial. Consider the reduced homology exact sequence of the pair $(K,M_D)$, where $K=\ov S\setminus (Tub(D)\cup Tub(\E))$: $$0\raa H_2(K,\Z)\raa H_2(K,M_D,\Z)\raa H_1(M_D,\Z)\raa H_1(K,\Z)\raa H_1(K,M_D,\Z)\raa 0.$$ By the Lefschetz duality (cf. \cite[3.43]{Hat}) $$H_i(K,M_D,\Z)\cong H^{4-i}(K,M,\Z)=H^{4-i}(S',Sing\ S',\Z),$$ which for $i>1$ implies that $$H_i(K,M_D,\Z)\cong H^{4-i}(S',\Z)\cong H_{3-i}(S',\Z)$$ by the universal coefficient formula. This gives an exact sequence: $$0\raa H_2(K,\Z)\raa H_1(S',\Z)\raa H_1(M_D,\Z)\raa H_1(K,\Z)\raa H_2(S',\Z)\raa 0.$$ Consider the reduced homology exact sequence of the pair $(K,M)$: $$0\raa H_2(K,\Z)\raa H_2(K,M,\Z)\raa H_1(M,\Z)\raa H_1(K,\Z)\raa H_1(K,M,\Z)\raa \wt H_0(M,\Z)\raa 0.$$ Since $H_i(K,M,\Z)\cong H_i(S',Sing\ S',\Z)$ and $H_1(S',Sing\ S',\Z)=H_1(S',\Z)\oplus \wt H_0(\Sing S',\Z)$ we get $$0\raa H_2(K,\Z)\raa H_2(S',\Z)\raa H_1(M,\Z)\raa H_1(K,\Z)\raa H_1(S',\Z)\raa 0.$$ Since $H_2(S',\Z)=0$ by (v), we get $H_2(K,\Z)=0$. Now $|H_1(M_D,\Z)|=|d(D)|$ and $|H_1(M,\Z)|=|d(\E)|$ by (i), so we obtain the result easily. \end{proof}

\bcor\label{cor3:homology} With the notation as above one has:\benum[(i)]

\item $b_1(S_0)=b_2(S_0)=b_1(\E),\ b_3(S_0)=q,\ b_4(S_0)=0$,

\item $\chi(S_0)=1-q$, $\chi(S)=\#\E+1-b_1(\E)$, $\chi(\ov S)=\#D+\#\E+2-2 b_1(\E)$,

\item $\Sigma_{S_0}=h+\nu-2$ and $\nu\leq 1$,

\item $S'$ is affine and $NS_\Q(S_0)=0$,

\item $d(D)<0$,

\item if $\pi_1(S')=\{1\}$ then $S'$ is contractible. \eenum\ecor

\begin{proof} Part (i) follows from \ref{prop3:homology}(i)-(ii). Then (ii) is a consequence of \ref{prop3:homology}(iii) and the equality $\chi(S_0)=\chi(S')-q=1-q$. By \ref{prop3:homology}(iv) $H_2(i_{D\cup\E})$ is surjective, so $NS_\Q(S_0)=0$ and then by \ref{cor2:Surf/Div} $S'$ is affine, which gives (iv). Since $H_2(i_{D\cup\E})$ is injective, the Hodge index theorem implies that the signature of $Q(D+\E)$ is $(1^-,\#(D+\E)^+)$, hence $$d(\E)d(D)=d(D+\E)<0,$$ which proves (v). For (iii) note that since $b_2(\ov S)=b_2(D\cup \E)$, Fujita's equation (sec. \ref{ssec:rulings}) yields $\Sigma_{S_0}=h+\nu-2$. If $\nu>1$ then the numerical equivalence of fibers of a $\PP^1$-ruling gives a numerical dependence of components of $D+\E$, hence $d(D+\E)=0$, a contradiction with (v). If $\pi_1(S')=\{1\}$ then by \ref{prop3:homology}(v) and the Hurewicz theorem all homotopy groups of $S'$ vanish, so Whitehead's theorem implies (vi). \end{proof}

\subsection{Birational type and logarithmicity}\label{ssec:bir type and logarithmicity}

By \cite[Theorem 1.1]{PS-rationality} it is known that singular $\Q$-homology planes which have at most quotient singularities are rational. We will see that this is not true in general without restrictions on the character of singularities. We describe the birational type of $S'$ and prove some general properties of the singular locus.

\bprop\label{prop3:general is logarithmic}Let $S_0$ be the smooth locus of a singular $\Q$-homology plane $S'$. \benum[(i)]

\item If $\ovk(S_0)=2$ then $S'$ is logarithmic and $\#\Sing S'=1$.

\item If $\ovk(S_0)=0$ or $1$ then either $\#\Sing S'=1$ or $\#\Sing S'=2$ and $\E_1=\E_2=[2]$,

\item If $\ovk(S_0)=-\8$ then either $S'\cong \C^2/G$ for some finite small non-cyclic subgroup $G<GL(2,\C)$ or $S'$ is affine-ruled and its singularities are cyclic.
\eenum\eprop

\begin{proof} (i)-(ii) Let $(S_m,D_m)$ be an almost minimal model of $(\ov S,D+\E)$. Since $S'$ is affine, an almost minimal model $S_m-D_m$ of $S_0$ is isomorphic to an open subset of $S_0$ satisfying $\chi(S_m-D_m)\leq \chi(S_0)=1-q$ (see \cite[2.8]{Palka-k(S_0)=0}). For each connected component of $D_m$ being a connected component of $\Bk D_m$ (hence contractible to quotient singularity) denote the local fundamental group of the respective singular point $P$ by $G_P$ and the set of such points by $Q$. By the Kobayashi inequality (see \cite{Langer} or \cite[2.5(ii)]{Palka-k(S_0)=0} for a generalization, which we use here) $$\frac{1}{3}((K_{S_m}+D_m)^+)^2\leq \chi(S_m-D_m)+\sum_{P\in Q}\frac{1}{|G_P|}\leq 1-q+\frac{\#Q}{2}\leq 1-\frac{q}{2}.$$ If $\ovk(S_0)=2$ then we get $q=1$ and $0<\ds\sum_{P\in Q}\frac{1}{|G_P|}$, so there is a unique singular point on $S'$ and it is of quotient type. If $\#\Sing S'>1$ and $\ovk(S_0)\geq 0$ then we get $q=2$ and $$1\leq \frac{1}{|G_{P_1}|}+\frac{1}{|G_{P_2}|},$$ so $|G_{P_1}|=|G_{P_2}|=2$.

(iii) If $S'$ is affine-ruled then it has only cyclic singularities by \cite{Miyan-cylinderlike}. If $S_0$ is not affine-ruled then by \cite{MiTs-PlatFibr} (see also \cite[2.5.1]{Miyan-OpenSurf}) it contains a Platonically $\C^*$-fibred open subset $U$. It is known that $U\cong (\C^2-0)/G$ for some small non-cyclic $G<GL(2,\C)$. Moreover, $S'\setminus U$ is a disjoint sum of affine lines, so since $$0\leq \chi(S')-\chi(U)=\chi(S')=1-\#\Sing S',$$ $S'$ has a unique singular point and $S_0=U$. Since $S'$ is normal, global regular functions on it are the same as the ones on $\C^2/G$. It follows that $S'\cong \C^2/G$, because $S'$ is affine. \end{proof}

\bprop\label{prop3:generalities on ovS-D-E} With the notation as above one has: \benum[(i)]

\item  $\ov S$ is $\PP^1$-ruled over a curve of genus $\frac{1}{2}b_1(D)=\frac{1}{2}b_1(\E)$,

\item if $\ovk(S')\geq 0$ then $S'$ is rational and has topologically rational singularities (cf. section \ref{sec:preliminaries}),

\item both $\E$ and $D$ are forests with at most one nonrational component,

\item if $\E$ consist only of $(-2)$-curves then $\ovk(S')=\ovk(S_0)$.

\eenum\eprop

\begin{proof}

(i)-(iii) We have $b_1(\E)=b_1(D)=b_1(\ov S)$ by \ref{prop3:homology}(iii), so if $b_1(\E)=0$ then we are done. We can therefore assume that $b_1(\E)\neq 0$. Suppose $\ovk(S)=-\8$. Then $S$ is affine-ruled (i.e. $\C^1$-ruled), because $D$ is connected, so we need only to prove (iii). Let $\ov S\to B$ be a $\PP^1$-ruling extending the affine ruling of $S$. Then $D$ is a tree and has a unique nonrational component as the horizontal section. Since $b_1(\E)\neq 0$, $\E$ has a horizontal component $E_0$. Clearly, $g(E_0)\geq g(B)$, so $b_1(E_0)\geq b_1(B)$. However, $b_1(B)=b_1(D)=b_1(\E)$, so $b_1(E_0)=b_1(\E)$ and $E_0$ is the unique horizontal component of $\E$, hence $\E$ is a forest. Thus we can assume that $\ovk(S_0)\geq \ovk(S)\geq 0$. Suppose $\ovk(S_0)=1$. Then, since $S'$ does not contain complete curves, by \cite[2.3]{Kawamata2} $S_0$ is $\C^*$-ruled and this ruling does not extend to $S'$ ($\E$ would be a rational forest then). Thus some resolution of $S'$ (not necessarily the minimal resolution $S$) is affine-ruled, which implies $\ovk(S)=-\8$, a contradiction. By \ref{prop3:general is logarithmic}(i) $\ovk(S_0)\neq 2$. Thus we are left with the case $\ovk(S)=\ovk(S_0)=0$. We argue as in the proof of \ref{prop3:homology}(iii)-(iv) that $b_1(\E)=0$, a contradiction.

(iv) We have to prove that $\ovk(S_0)\leq \ovk(S)$. If $\E$ consists of $(-2)$-curves then $(K+D)\cdot E_i=0$ for each irreducible component $E_i$ of $\E$. If $T$ is an effective divisor linearly equivalent to $n(K+D+\E)$ then, since $Q(\E)$ is negative definite, $T-n\E$ is effective by \ref{lem2:Eff-NegDef is Eff} and we are done. \end{proof}

Note that \ref{cor3:homology}(iv) and \ref{prop3:generalities on ovS-D-E}(i) establish \ref{thm:main result1}(1). Let us denote by \ref{thm:main result1}(2')\label{thm:main result1prim} a weaker version of the theorem \ref{thm:main result1}(2) with parts (a) and (b) replaced by the following weaker statements:
\begin{enumerate} \item[(a')] It has a $\C^1$- or a $\C^*$-ruling. \item[(b')] Its smooth locus has a $\C^*$-ruling which does not extend to a ruling of the $\Q$-homology plane. \end{enumerate} We will prove the original version after proving \ref{thm:main result2}.

\begin{proof}[Proof of \ref{thm:main result1}(2')]

By \ref{prop3:general is logarithmic}(iii) we may assume $\ovk(S_0)\geq 0$. By \cite[1.4]{Palka-k(S_0)=0} we have only to consider the case when singularities of $S'$ are not topologically rational, i.e. we may assume $b_1(\E)\neq 0$. Then $\E$ is connected by \ref{prop3:general is logarithmic} and $b_1(D)\neq 0$ by \ref{prop3:homology}(iii). Thus by \ref{prop3:generalities on ovS-D-E} $D$ and $\E$ are trees with unique non-rational components. Let $(\wt S,\wt D+\wt E)$ be an almost minimal model of $(\ov S,D+\E)$. By the affiness of $S'$, the $(-1)$-curves contracted in the process of minimalization intersect the image of $D$. Since $D$ is not negative definite, its image cannot be an admissible chain, hence by \ref{rem:almost minimalization} the $(-1)$-curves do not intersect $\E$. By \cite[8.8]{Fujita} $\wt D$ and $\wt E$ are disjoint smooth elliptic curves. By \ref{prop3:generalities on ovS-D-E}(i) $\ov S$ is $\PP^1$-ruled over a smooth elliptic curve, so L\"{u}roth's theorem implies that every rational curve in $\ov S$ is vertical. In particular, $(-1)$-curves contracted in the process of minimalization are vertical, so the ruling descends to $\wt S$ and the number of horizontal components of $D+\E$ and $\wt D+\wt E$ is the same. For a general fiber $f$ we get $$-2+f\cdot (D+\E)=f\cdot  K_{\wt S}+ f\cdot \wt D+f\cdot \wt E=f\cdot \Bk(\wt D+\wt E)=0,$$ because all components contained in $\Supp\Bk(\wt D+\wt E)$ are rational, hence vertical. Thus $f\cdot (D+\E)=2$, so $S_0$ is $\C^*$-ruled and we are in case (b'). \end{proof}

\section{\label{sec:nonlogarithmic S'} Non-logarithmic $S'$}

\subsection{Existence of a non-extendable $\C^*$-ruling} The following result strengthens the proposition \ref{prop3:generalities on ovS-D-E}(ii) and is the main step towards \ref{thm:main result2}.

\bthm\label{thm3:nonlogarithmic are special} If a singular $\Q$-homology plane is not logarithmic then its smooth locus has a unique $\C^*$-ruling. This $\C^*$-ruling does not extend to a ruling of the plane. Moreover, the Kodaira dimension of the plane is negative and the Kodaira dimension of its smooth locus is zero or one. \ethm

\begin{proof} By \ref{prop3:general is logarithmic} $\ovk(S_0)\in \{0,1\}$ and $\E$ is connected. By \ref{thm:main result1}(2') (see the remarks after \ref{prop3:generalities on ovS-D-E}) we may assume that $S_0$ is $\C^*$-ruled. We will first show that this ruling cannot be extended to a ruling of $S'$. Consider a minimal completion $(\ov S,D+\E,\pi)$ of a $\C^*$-ruling of $S_0$. It is enough to show that $\E_h\neq 0$. Suppose $\E_h=0$. Then $D_h$ consists either of two 1-sections or of one 2-section. In particular, it can intersect only those components which have multiplicity one or two. In the second case $\#D_h=1$ and the point of intersection is a branching point of $\pi_{|D_h}$. The exceptional divisor $\E$ is vertical, so $\ov S$ and $D$ are rational by \ref{prop3:generalities on ovS-D-E}(i). Let $F$ be the singular fiber containing $\E$ and let $D_v$ be the divisor of $D$-components of $F$. By \ref{cor3:homology}(iii) we have $\nu\leq 1$ and $\Sigma_{S_0}=\#D_h+\nu-2\leq 1$, so $\sigma>1$ for at most one fiber of $\pi$. We obtain successive restrictions on $F$ eventually leading to a contradiction. We use \ref{lem2:fiber with one exc comp} repeatedly.

\bcl The $(-1)$-curves of $F$ are $S_0$-components.\ecl

Suppose $F$ contains a $(-1)$-curve $D_0\subseteq D$. Since $\pi$ is minimal, the divisor $D_h$ intersects $D_0$, so either $\mu(D_0)=1$ or $\mu(D_0)=2$. Moreover, $D_0$ can be a tip of $F$ only if $D_h$ intersects it in two distinct points. In particular, we see that $D_v$ contains a component of multiplicity one and does not contain more $(-1)$-curves. We have $\Sigma_{S_0}=0$. Indeed, if $\Sigma_{S_0}=1$ then $\#D_h=2$ and $\nu=1$, so by simply connectedness of $D$ at most one horizontal component of $D$ meets $D_0$. However, in this case $\mu(D_0)=1$, so $D_0$ is a tip of $F$, a contradiction. The unique $S_0$-component $C$ of $F$ is exceptional, otherwise $D_0$ would be the unique $(-1)$-curve of $F$, which would imply that $F=[2,1,2]$ with no place for $\E$. Clearly, there are no more $(-1)$-curves in $F$. Let us make a connected sequence of blowdowns starting from $D_0$ until the number of $(-1)$-curves decreases. Since $\E\cap D=\emptyset$, in this process we do not touch $C+\E$ (first we would touch $C$, and then $C$ becomes a $0$-curve). Let $F'$ be the image of $F$, we can write $\un F'-C=D'+\E$, where $D'$ is the image of $D_v$. Since $C+\E$ is not touched, $D'\neq 0$. We know that $D_v$ contains a component of multiplicity one, so the same is true for $D'$. It follows that $\E$ is a chain, a contradiction.

\bcl $F$ contains two $(-1)$-curves. \ecl

Suppose $F$ has a unique $(-1)$-curve $C$. It still might contain two $S_0$-components. Write $\un F-C=A+B$, where $A$ and $B$ are disjoint, connected, and $B$ is a chain (possibly empty). By our assumption on $\E$ we have $\E\subseteq A$ and $A$ is not a chain. Thus $B$ contains only $S_0$- and $D$-components. Note that by affiness of $S'$ each $S_0$-component meets $D$. By connectedness of $D$ either $B\cdot D_h>0$ or $B=0$. Suppose $B\neq 0$. Then $B$ contains a curve with $\mu\leq 2$, so $F$ has two branches, the first equal to $[2,k,2]$ for some $k>1$. In this case $\E$ is an admissible fork of type $(2,2,n)$, a contradiction. Thus $B=0$. Since $\E$ is not an admissible fork, $\mu(C)>2$. If follows that $D_h\cdot C=0$ and $C$ intersects $D$ in one point, belonging to some component $D_1\subseteq D$. Since $D$ is connected, there is a chain $T\subseteq D_v$ containing $D_1$ and some $D_2\subseteq D$ intersecting $D_h$. We have $\mu(D_2)\leq 2$ and since $\E$ contains a branching component of $F$, $D_2$ does not belong to the first branch of $F$. In fact it follows that $D_2$ belongs to the second branch and $\E$ is an admissible fork of type $(2,2,n)$, a contradiction.

\bcl Both $(-1)$-curves of $F$ intersect $\E$.\ecl

Let $C_1$ and $C_2$ be the $(-1)$-curves of $F$. They are unique $S_0$-components of $F$, because $\sigma(F)\leq 2$. Now $D_h$ consists of two $1$-sections, which can intersect $F$ only in components of multiplicity one. Suppose one of $C_i$'s, say $C_2$, does not meet $\E$. Then $D_v\neq 0$, because $C_2$ meets some component of $F$. As in the proof of Claim $1$ we make a connected sequence of blowdowns starting from $C_2$ until there is only one $(-1)$-curve left, we denote the image of $F$ by $F'$. Again in this process we do not touch $C_1+\E$, so we can write $\un F'-C_1=D'+\E$, where $D'$ is the image of $D_v$. Since $D'$ intersects the image of $D_h$, it contains a component of multiplicity one. It follows that $\E$ is a chain, a contradiction.

\bcl There are no $D$-components in $F$.\ecl

We can write $\un F-C_1-C_2=\E+D'+D''$, where $D_v=D'+D''$, $D'$ and $D''$ are connected and disjoint. Suppose $D'\neq 0$. One of $C_i$'s, say $C_1$, meets $D'$. Contract $C_2$ and subsequent $(-1)$-curves until the number of $(-1)$-curves decreases. Clearly, $C_1+D'$ is not affected in this process. Denote the image of $F$ by $F'$ and let $U$ be the image of $D''+C_2+\E$. Now $F'$ is a fiber with a unique $(-1)$-curve and since both $C_2+D''$ and $C_1+D'$ intersect $D_h$, we infer that both $U$ and $C_1+D'$ contain components of multiplicity one. Thus $F'$ is a chain. Consider the reverse sequence of blowups recovering $F$ from $F'$. The fiber $F$ is not a chain, so a branching curve is produced. It follows that $D''+C_2$ contains no curve of multiplicity one, so $D_h\cdot (D''+C_2)=0$, a contradiction.

\medskip The last claim implies that $D_h$ intersects both $C_i$'s, so they have multiplicity one, hence are tips of $F$. It follows that $F$ is a chain. Thus $\E$ is a chain, a contradiction. This finishes the proof that no $\C^*$-ruling of $S_0$ can be extended to a ruling of $S'$. We see also that $S$ is affine-ruled, hence $\ovk(S')=-\8$.

Note that since $\E_h\neq 0$, we have $\nu=0$, so by \ref{cor3:homology}(iii) $\Sigma_{S_0}=0$, i.e. every fiber of $\pi$ has exactly one $S_0$-component. By \cite[7.6]{Fujita} every singular fiber is columnar, so $\E_h$ and $D_h$ are the unique branching components of $\E$ and $D$ respectively ($\E_h$ is branching, as $\E$ is not a chain). All non-branching components of $D+\E$ have negative self-intersection, so $(\ov S,D+\E)$ is actually the unique snc-minimal completion of $S_0$. Indeed, if $(X,W)$ with birational morphisms $\ov \eta\:(X,W)\to (\ov S,D+\E)$ and $\wt \eta\:(X,W)\to (\wt S,\wt D+\wt E)$ is a minimal correspondence between $(\ov S,D+\E)$ and some other snc-minimal completion of $S_0$ then $\wt \eta$, and hence also $\ov \eta$, is an identity, as any $(-1)$-curve $L\subseteq W$ contracted by $\wt \eta$ would give a non-branching component $\ov \eta_*L\subseteq D+\E$ with non-negative self-intersection.

Suppose now that $S_0$ has a second $\C^*$-ruling which cannot be extended to a ruling of $S'$. By the remarks above the ruling extends to a $\PP^1$-ruling $\pi'$ of $\ov S$ and the components $D_h$, $\E_h$, which are by definition horizontal for $\pi$, are also horizontal for $\pi'$. But then $D-D_h+\E-\E_h$ is vertical for both rulings. Thus, if $F$ and $F'$ are fibers of $\pi$ and $\pi'$ respectively then $F-F'$ intersects trivially with all components of $D+\E$. By \ref{cor3:homology}(iv) $F\equiv F'$, so $F\cdot F'=F^2=0$, i.e. $\pi=\pi'$. \end{proof}

\subsection{Construction and properties}\label{ssec:nonlogarithmic-construction} We now proceed to prove the remaining parts of \ref{thm:main result2}. From now on $S'$ is a singular $\Q$-homology plane with a $\C^*$-ruled smooth locus $S_0$. Note that $\ovk(S_0)\neq 2$ by the easy addition theorem \cite[10.4]{Iitaka}. We also assume that the ruling is \emph{non-extendable}, meaning that it does not extend to a ruling of $S'$. By \ref{thm3:nonlogarithmic are special} this is the case if $S'$ is non-logarithmic. Let $(\ov S,D+\E,p)$ be a minimal completion of such a $\C^*$-ruling of $S_0$, where $\E$ is an exceptional locus of some resolution of singularities of $S'$. We have $D_h\neq 0$, otherwise $D$ would be vertical, which contradicts the affiness of $S'$. Since $p_{|S_0}$ does not extend to a ruling of $S'$, we have $\E_h\neq 0$. So $p$ is untwisted and $S$ is affine-ruled, which gives $\ovk(S')=-\8$. Let $N=-\E_h^2$ and let $F_1,F_2,\ldots,F_n$ be all the columnar fibers of $p$. Let $E_i\subseteq F_i$ be connected components of $\E-\E_h$. Let $C_i$ be the unique $(-1)$-curve of $F_i$, put $\mu_i=\mu(C_i)$. Note that $\mu_i$ is the denominator of the reduced form of the fraction $\wt e(E_i)$ (cf. remark after \ref{def2:(un)twisted and columnar fiber}). Denote the base of $p$ by $B$. By \ref{prop3:generalities on ovS-D-E}(i) the rationality of one of $\ov S$, $\E$, $D$ or $B$ implies the rationality of all others.

\begin{figure}[h]\centering\includegraphics[scale=0.4]{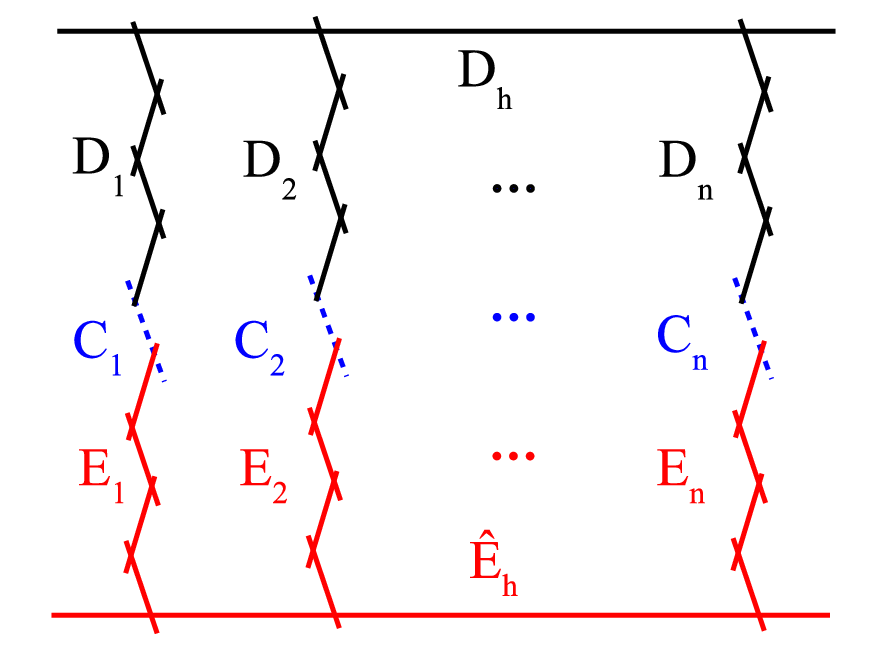}\caption{Non-extendable $\C^*$-ruling}  \label{fig:nonextendable}\end{figure}

\blem\label{lem5:nonextendable ruling properties} Singular fibers of $p$ are columnar and $\sum_{i=1}^n \wt e(E_i)<N$ (see Fig.~\ref{fig:nonextendable}). There exists a line bundle $\mathcal{L}$ over $B$ with $\deg \mathcal{L}=-N<0$ and a proper birational morphism $\ov S\to\PP(\mathcal{O}_B \oplus\mathcal{L})$, such that $p$ is induced by the projection of $\PP(\mathcal{O}_B \oplus\mathcal{L})$ onto $B$.\elem

\begin{proof} Since $\E\cap D=\emptyset$, we have $\nu=0$. By \ref{cor3:homology}(iii) $\Sigma_{S_0}=0$, so every fiber has exactly one $S_0$-component. By \cite[7.6]{Fujita} every singular fiber is columnar. We contract all singular fibers to smooth fibers (i.e. we contract subsequently their $(-1)$-curves) without touching $\E_h$. Denote the image of $\ov S$ by $\wt{S}$ and the image of $D_h$ by $\wt{D}_h$. Then $\E_h$ is disjoint from $\wt{D}_h$. Since $\E_h^2=-N<0$, we can write $\wt{S}=\PP(\mathcal{O}_B\oplus \mathcal{L})$ for a line bundle $\mathcal{L}$ with $deg\ \mathcal{L}=-N<0$ (see \cite[V.2]{Har}). Now $\E_h$ and $\wt D_h$ are sections coming from the direct summands of the bundle. The matrix $Q(\E)$ is negative definite, so (cf. \cite[2.1.1]{KR-ContrSurf}) $$0<\det Q(-\E)=d(E_1)d(E_2)\ldots d(E_n)(-\E_h^2-\sum_{i=1}^n \wt e(E_i)),$$ hence $\sum_{i=1}^n \wt e(E_i)<N$. \end{proof}

\bcor\label{cor5:nonextendable is contractible} $S'$ is contractible. \ecor

\begin{proof} By \ref{lem5:nonextendable ruling properties} singular fibers of $p$ are columnar, so in each fiber there is a component of $\E$ of multiplicity one, hence by \cite[5.9, 4.19]{Fujita} the embedding $\E_h\to S$ induces an isomorphism $\pi_1(\E_h)\to \pi_1(S)$. Thus by \ref{prop3:homology}(v)-(vi) and Whitehead's theorem $S'$ is contractible. \end{proof}

\bcon\label{con5:nonextendable} Pick $n\in \N$ and for each $i=1,\ldots,n$ choose a number $\wt e_i\in \Q\cap (0,1)$. Choose a positive integer $N$, such that $\sum_{i=1}^n \wt e_i<N$. Let $B$ be a complete curve of genus $g(B)$, such that $g(B)>0$ if $n$ was chosen smaller than $3$. Define $\wt S=\PP(\cal O_B\oplus \cal L)$, where $\cal L$ is a line bundle over $B$ of degree $\deg \cal L=-N$. Let $\wt p\:\wt S\to B$ be the induced $\PP^1$-fibration. Denote the sections induced by inclusions of the direct summands $\mathcal O_B$ and $\cal L$ by $\E_h$ and $\wt D_h$. Then $\E_h^2=-N$ and $\wt{D}_h^2=N$. Choose $n$ distinct points $x_1,\ldots,x_n\in\wt D_h$ and blow up each point once. For each $i$ make a connected sequence of subdivisional blowups creating over $\wt p(x_i)$ a columnar fiber $F_i$ with $\wt e(E_i)=\wt e_i$. Denote the birational transform of $\wt D_h$ by $D_h$. Write $\un {F_i}=E_i+C_i+D_i$ where $C_i^2=-1$, $E_i$ and $D_i$ are connected chains and $D_i\cap \E_h=\emptyset$. Let $\mu_i$ be the multiplicity of $C_i$ in $F_i$. Fix a natural order on each $E_i$ and $D_i$ treated as twigs of $\E=E_1+\ldots+E_n+\E_h$ and $D=D_1+\ldots+D_n+D_h$ respectively. Denote the obtained surface by $\ov S$ and the induced $\PP^1$-ruling by $p$. Define $S=\ov S-D$, $S_0=S-\E$ and $S'=S/\E$ (as a topological space). We will show below that $NS_\Q(S_0)=0$, hence by \ref{cor2:Surf/Div} $S'$ and the quotient morphism can be realized in the algebraic category.\econ

\bsrem The additional assumption that $g(B)>0$ if $n<3$ is justified as follows. If $g(B)=0$ and $n<3$, then $D$ is a non-admissible chain, so by blowing up and down on $D$ we may assume it contains a $0$-curve as a tip. The tip induces an affine ruling of $S_0$, hence $\ovk(S_0)=-\8$. Moreover, $\E$ is a chain so the singularity is cyclic. $\Q$-homology planes of this type are described in \cite[2.7, 2.8]{MiSu-hPlanes}.\esrem

\bthm\label{thm5:nonextendable} The surface $S'$ constructed in \ref{con5:nonextendable} is a contractible surface of negative Kodaira dimension. Moreover, each non-logarithmic $\Q$-homology plane (or more generally, each non-affine-ruled $\Q$-homology plane whose smooth locus admits a non-extendable $\C^*$-ruling) can be obtained by construction \ref{con5:nonextendable}. The Kodaira dimension of the smooth locus is determined by the sign of the number $$\alpha=n-2+2g(B)-\sum_{i=1}^n\frac{1}{\mu_i}$$ (i.e. $\ovk(S_0)=-\8$ for $\alpha<0$, $0$ for $\alpha=0$ and $1$ for $\alpha>0$). The snc-minimal completion and the pair $(B,\cal L)$ used in the construction are determined uniquely by the isomorphism type of $S'$.\ethm

\begin{proof} The assumption that $S'$ is not affine-ruled excludes the case when $g(B)=0$ and $n\leq 2$, as was done in the construction. It follows from \ref{lem5:nonextendable ruling properties} and \ref{thm3:nonlogarithmic are special} that  if $S'$ is not affine-ruled but admits a non-extendable $\C^*$-ruling then it can be obtained by construction \ref{con5:nonextendable}. The matrix $Q(\E-\E_h)$ is negative definite and $$d(\E)=d(E_1)d(E_2)\cdots d(E_n) (N-\sum_{i=1}^n \wt e_i)>0,$$ so by Sylvester's theorem $Q(\E)$ is negative definite. We have $$d(D)=d(D_1)d(D_2)\cdots d(D_n) (-N+n-\sum_{i=1}^n (1-\wt e_i))=-d(\E)$$ by the remark after \ref{def2:(un)twisted and columnar fiber}, so $d(D)\neq 0$. It follows that the classes of irreducible components of $D+\E$ are independent in $NS_\Q(\ov S)$, hence are a basis because $b_2(\ov S)=\#D+\#\E$. We apply \ref{cor2:Surf/Div} and infer that $S'$ is normal and affine. By Iitaka's easy addition theorem $\ovk(S_0)\leq 1$. The divisor $K_{\ov S}+D+\E+\ds\sum_{i=1}^n C_i$ intersects trivially with all vertical components, so it is numerically equivalent to a multiple of a general fiber $f$. Intersecting with $D_h$ we get $$K_{\ov S}+D+\E+\sum_{i=1}^n C_i\equiv (2g(B)-2+n) f.$$ Putting $G_i=\frac{1}{\mu_i}F_i-C_i$ we get $K_{\ov S}+D+\E\equiv \alpha f+\ds\sum_{i=1}^n G_i$. Since $\ds\sum_{i=1}^n G_i$ is effective, vertical and has a negative definite intersection matrix, by \ref{lem2:Eff-NegDef is Eff} we get $\kappa(K_{\ov S}+D+\E)=\kappa(\alpha f)$, so $\ovk(S_0)$ is determined by the sign of $\alpha$ as stated.

Now we check that $S'$ is $\Q$-acyclic (then it is contractible by \ref{cor5:nonextendable is contractible}). We know from the above that the map $H_2(D+\E)\to H_2(\ov S)$ induced by inclusion is an isomorphism. Clearly, $H_1(D)\to H_1(\ov S)$ and $H_1(\E)\to H_1(\ov S)$ are monomorphisms, because they are monomorphisms after composing with $H_1(p)$. The exact sequence of the pair $(D,\ov S)$ gives $b_4(S)=b_3(S)=0$, $b_2(S)=\#\E$ and $b_1(S)=b_1(\ov S)=b_1(B)$. Then the exact sequence of the pair $(\E,S)$ gives $$b_1(S')=b_2(S')=b_3(S')=b_4(S')=0.$$ Since we assumed that $g(B)>0$ if $n<3$, $S'$ is singular.

All non-branching rational curves contained in $D$ have negative self-intersection, so the smooth completion of $S$ is unique up to isomorphism (it is snc-minimal, as $B$ is branching if $g(B)=0$). Suppose $S_1'\cong S_2'$ are two surfaces constructed as in \ref{con5:nonextendable}, we will use indices $1,2$ consequently to distinguish between objects occurring in the intermediate steps of the construction. Since all rational non-branching components of $D+\E$ have negative self-intersection, the isomorphism extends to an isomorphism of completions $\Phi\:(\ov S_1,D+\E)\to (\ov S_2,D+\E)$. Now the argument from the proof of \ref{thm3:nonlogarithmic are special} shows that there is at most one non-extendable $\C^*$-ruling of $S_0$, so up to composition with an automorphism of $\ov S$ induced by an automorphism of $B$ we can assume that $\Phi$ preserves fibers , so in particular it fixes all components of $\E+D$. Then $\Phi$ induces an isomorphism of $B$-schemes $\wt S_1$ and $\wt S_2$. Thus $$\cal O_B\oplus \cal L_1\cong (\cal O_B\oplus \cal L_2)\otimes \cal E$$ for some line bundle $\cal E$ of degree zero. It follows that $\deg \cal L_2\otimes\cal E<0$, so non-vanishing constant sections of $\cal O_B$ (on the left hand side) are sections of $\cal E$, which gives $\cal E\cong \cal O_B$. Thus $$\cal O_B\oplus\cal L_1\cong \cal O_B\oplus\cal L_2$$ which after taking second exterior power gives $\cal L_1\cong\cal L_2$. \end{proof}

\brem Note that the non-affine ruled $\Q$-homology planes $S'$ with $\ovk(S_0)=-\8$ can also obtained by the construction above. Indeed, the smooth locus of such $S'$ is Platonically $\C^*$-fibred (see the proof of \ref{prop3:general is logarithmic}(iii)) and this fibration cannot be extended to a ruling of $S'$. \erem

\bcor\label{cor5:nonextendable-singularities} Let $P\in S'$ be the unique singular point of a singular $\Q$-homology plane $S'$, whose smooth locus has a non-extendable $\C^*$-ruling. Then with the notation as above: \benum[(i)]

\item $P$ is a topologically rational singularity if and only if $B\cong\PP^1$,

\item if $\ovk(S_0)=-\8$ then $g(B)=0$, $n\leq 3$ and $S'$ is logarithmic. If additionally $n>2$ (as assumed in the construction) then $(\mu_1,\mu_2,\mu_3)$ is up to order one of the Platonic triples (i.e. triples $(x_1,x_2,x_3)$ of positive integers satisfying $\sum_{i=1}^3\frac{1}{x_i}>1$), so $S_0$ has a structure of a Platonic fibration.

\item if $\ovk(S_0)\geq 0$ then $S'$ is not logarithmic,

\item $\ovk(S_0)=0$ if and only if either \benum[(a)] \item $g(B)=1$ and $n=0$ or \item $g(B)=0$, $n=4$ and $\mu_1=\mu_2=\mu_3=\mu_4=2$ or \item $g(B)=0$, $n=3$ and $(\mu_1,\mu_2,\mu_3)$ is up to order one of $(2,3,6),(2,4,4),(3,3,3)$.\eenum \eenum\ecor

\begin{proof}
(ii) If $\alpha<0$ then $$\frac{n}{2}\leq\sum_{i=1}^n(1-\frac{1}{\mu_i})<2(1-g(B)),$$ so $g(B)=0$ and $n\leq 3$. Suppose $n=3$. Then $\sum_{i=1}^3\frac{1}{\mu_i}>1$, so $(\mu_1,\mu_2,\mu_3)$ is up to order one of the Platonic triples.

(iii) If $S'$ is logarithmic then $\E$ is either a chain or an admissible fork. In the first case $n\leq 2$ and in the second $n=3$ and $\sum_{i=1}^3 \frac{1}{\mu_i}>1$. In both cases $\alpha<0$, so $\ovk(S_0)=-\8$.

(iv) Assume $\alpha=0$. For $n=0$ we get $g(B)=1$. Assume $n>0$. We have $$\frac{n}{2}\leq\sum_{i=1}^n(1-\frac{1}{\mu_i})=2(1-g(B)),$$ so we get $g(B)=0$ and $n\in \{3,4\}$. We have then $\sum_{i=1}^3\frac{1}{\mu_i}=1$ if $n=3$ and $\sum_{i=1}^4\frac{1}{\mu_i}=2$ if $n=4$, which gives (b) and (c). Conversely, in each case $\alpha=0$. \end{proof}

\bex\label{ex4:nonextendable-singularities} Suppose $n\geq 3$, $N\geq 1$, $\wt e_1,\ldots,\wt e_n\in\Q\cap (0,1)$, $\wt e_1+\ldots+\wt e_n<N$ and $B\cong \PP^1$. Let $P$ be the unique singular point of $S'$ constructed as in \ref{con5:nonextendable}.\benum[(i)]

\item If $N\geq n$ then $P\in S'$ is a rational singularity.

\item If $N<n-1$ then $P\in S'$ is a topologically rational but not a rational singularity. \eenum \eex

\begin{proof}

(i) We have $\sum_{i=1}^n \wt e_i<n\leq N$. The fundamental cycle $Z_f$ of $\E$, which is the smallest nonzero effective divisor $Z_f\subseteq \E$, such that $Z_f\cdot E'\leq 0$ for each irreducible $E'\subseteq \E$, equals $\E$ in this case. Then $p_a(Z_f)=0$, so $P$ is a rational singularity by \cite[Theorem 3]{Artin-Rational}.

(ii) Let $Z=\E+\beta \E_h$, where $\beta=\lceil \frac{n}{N}\rceil-1$ ($\lceil x \rceil$ is defined as the smallest integer not smaller than $x$). Then $$p_a(Z)=\beta(n-1-\frac{1}{2}(\beta+1)N),$$ which is non-negative. Indeed, $\beta\geq 1$ and $$(\beta+1)N<(\frac{n}{N}+1)N=n+N\leq 2n-1,$$ so $p_a(Z)\geq 0$. However, if $p_a(Z)=0$ then the equality $\lceil\frac{n}{N}\rceil\cdot N=2n-2$ gives $n<N+2$, so $n=N+1$, a contradiction. It follows from \cite[Proposition 1]{Artin-Rational} that $P$ is not a rational singularity. \end{proof}

\bsrem As for (ii) note that for instance if $n>N+1$ and $E_i=[x_i]$ with $x_i\geq \frac{n}{N}$, not all equal $\frac{n}{N}$ then the condition $\wt e_1+\ldots+\wt e_n<N$ is satisfied. In general the fundamental cycle can be computed using \cite[Proposition 4.1]{Laufer}. \esrem

Note that $\wt S=\PP(\cal L_1\oplus\cal L_2)$ admits a natural $\C^*$-action fixing precisely the sections coming from the inclusion of direct summands. (Each $v\in \wt S$ can be written as $v=[u_1+u_2]$, where $u_i\in\cal L_{|\pi(v)}$ and $\pi$ is the projection onto the base, and then the action can be written as $t*[u_1+u_2]=[u_1+tu_2]$.) This action lifts to a $\C^*$-action on $\ov S$ constructed in \ref{con5:nonextendable}, because centers of successive blowups creating $\ov S$ belong to fixed loci of successive liftings. By \cite[1.1]{Pinkham-C*-surfaces} this gives the following corollary, which together with \ref{thm3:nonlogarithmic are special} and \ref{thm5:nonextendable} establishes \ref{thm:main result2}.

\bcor\label{cor4:nonlogarithmic are cones} A singular $\Q$-homology plane whose smooth locus admits a non-extendable $\C^*$-ruling (in particular any non-logarithmic) is a quotient of an affine cone over a smooth projective curve by an action of a finite group acting freely off the vertex of the cone and respecting the set of lines through the vertex.\ecor

We note that it follows from 5.8 loc. cit. that the unique singular point of $S'$ as in \ref{cor5:nonextendable-singularities} is a rational singularity if and only if $B$ is rational and $N>\frac{1}{k}(\sum_{i=1}^n\lceil k\wt e_i \rceil-2)$ for every natural number $k>0$.

\smallskip

\begin{proof}[Proof of \ref{thm:main result1}(2)] At the end of section \ref{ssec:bir type and logarithmicity} we proved 1.1(2'). As usually, let $S'$ be a singular $\Q$-homology plane with smooth locus $S_0$ not of general type. By \ref{prop3:general is logarithmic}(iii) we may assume $\ovk(S_0)\geq 0$. In particular $S'$ is not $\C^1$-ruled. If $S'$ is $\C^*$-ruled then this ruling is an extension of a $\C^*$-ruling of $S_0$, so by \ref{thm:main result2}(ii) $S'$ is logarithmic, hence satisfies \ref{thm:main result1}(2)(a). Finally, if $S'$ is as in (b'), i.e. $S_0$ has a $\C^*$-ruling which does not extend to a ruling of $S'$, then by \ref{cor5:nonextendable-singularities}(iii) $S'$ is non-logarithmic, hence satisfies \ref{thm:main result1}(2)(b). \end{proof}

\bibliographystyle{amsalpha}
\bibliography{bibl}

\end{document}